\documentclass[a4paper,12pt]{amsart}

\usepackage{amsfonts}
\usepackage{amsmath}
\usepackage{amssymb}

\numberwithin{equation}{section}
\theoremstyle{plain}
\newtheorem{thm}{Theorem}[section]
\newtheorem{proposition}[thm]{Proposition}

\theoremstyle{definition}

\newtheorem{remark}[thm]{Remark}
\newtheorem{example}[thm]{Example}

\def\R{{\mathbb R}}
\def\C{{\mathbb C}}
\def\Rn{{{\mathbb R}^n}}

\def\Sph{{\mathbb S}}
\def\FT{{\mathcal F}}
\def\L2tx{{L^2(\R_t\times\R^n_x)}}

\def\p#1{{\left({#1}\right)}}
\def\b#1{{\left\{{#1}\right\}}}

\def\n#1{{\left\|{#1}\right\|}}
\def\abs#1{{\left|{#1}\right|}}
\def\jp#1{{\left\langle{#1}\right\rangle}}
\def\Im{\operatorname{Im}}

\def\supp{\operatorname{supp}}

\def\rank{\operatorname{rank}}

\def\ka{\kappa}

\begin{document}

\title{Recent progress in smoothing estimates for
evolution equations}

\author[Michael Ruzhansky]{Michael Ruzhansky}
\address{
  Michael Ruzhansky:
  \endgraf
  Department of Mathematics
  \endgraf
  Imperial College London
  \endgraf
  180 Queen's Gate, London SW7 2AZ
  \endgraf
  United Kingdom
  \endgraf
  {\it E-mail address} {\rm m.ruzhansky@imperial.ac.uk}
  }
\author[Mitsuru Sugimoto]{Mitsuru Sugimoto}
\address{
  Mitsuru Sugimoto:
  \endgraf
  Graduate School of Mathematics
  \endgraf
  Imperial College London
  \endgraf
  Nagoya 464-8602
  \endgraf
  Japan
  \endgraf
  {\it E-mail address} {\rm sugimoto@math.nagoya-u.ac.jp}
  }

\thanks{
The authors were supported by the Daiwa Anglo-Japanese Foundation.
The first author was also supported by the EPSRC Leadership Fellowship
EP/G007233/1. }

%
%


\begin{abstract}
This paper is a survey article of results and arguments
from authors' papers
\cite{RS6}, \cite{RS2}, and \cite{RS7},
and describes a new approach to global smoothing problems
for dispersive and non-dispersive evolution equations 
based on ideas of comparison principle and canonical transforms.
For operators $a(D_x)$ of order $m$
satisfying the dispersiveness condition $\nabla a(\xi)\neq0$,
the smoothing estimate
\begin{equation*}
\n{\jp{x}^{-s}|D_x|^{(m-1)/2}e^{ita(D_x)}
\varphi(x)}_{L^2\p{\R_t\times\R^n_x}}
\leq C\n{\varphi}_{L^2\p{\R^n_x}} 
\quad {\rm(}s>1/2{\rm)}
\end{equation*}
is established,
while it is known to fail for general non-dispersive operators.
Especially, time-global smoothing estimates for the operator
$a(D_x)$ with lower order terms are the benefit
of our new method.
For the case when the dispersiveness breaks, we suggest a form
\begin{equation*}
\n{\jp{x}^{-s}|\nabla a(D_x)|^{1/2}
e^{it a(D_x)}\varphi(x)}_{L^2\p{\R_t\times\R^n_x}}
\leq C\n{\varphi}_{L^2\p{\R^n_x}} \quad{\rm(}s>1/2{\rm)}
\end{equation*}
which is equivalent to the usual estimate in the dispersive case and
is also invariant under canonical transformations
for the operator $a(D_x)$.
It does continue to hold for a variety of non-dispersive
operators $a(D_x)$, where $\nabla a(\xi)$ may become zero on some set.
It is remarkable that our method allows us to carry out a global
microlocal reduction of equations to the translation invariance 
property of the Lebesgue
measure.
\end{abstract}

\maketitle

\section{Introduction}
This survey article is a collection of results and arguments
from authors' papers
\cite{RS6}, \cite{RS2}, and \cite{RS7}.

Let us consider the following Cauchy problem to the Schr\"odinger equation:
\[
\left\{
\begin{aligned}
\p{i\partial_t+\Delta_x}\,u(t,x)&=0\quad\text{in $\R_t\times\R^n_x$},\\
u(0,x)&=\varphi(x)\quad\text{in $\R^n_x$}.
\end{aligned}
\right.
\]
By Plancherel's theorem, the solution
$u(t,x)=e^{it\triangle_x}\varphi(x)$ preserves the $L^2$-norm of the initial
data $\varphi$, that is, we have
$\n{u(t,\cdot)}_{L^2(\R^n_x)}=\n{\varphi}_{L^2(\R^n)}$
for any fixed time $t\in\R$.
But if we integrate the solution in $t$, we get an extra gain of regularity
of order
$1/2$ in $x$.
For example we have the estimate
\[
\n{\langle x\rangle^{-s}|D_x|^{1/2} e^{it\Delta_x}\varphi}
_{L^2\p{\R_t\times\R^n_x}}
\le C\n{\varphi}_{L^2\p{\R^n_x}} \qquad (s>1/2)
\]
for $u=e^{it\Delta_x}\varphi$, where $\jp{x}=\sqrt{1+|x|^2}$,
and (a sharper version of) this estimate was first given by
Kenig, Ponce and Vega \cite{KPV1}.
This type of estimate is called a smoothing estimate,
and its local version was first
proved by Sj\"olin \cite{Sj}, Constantin and Saut \cite{CS}, and Vega \cite{V}.
We remark that, historically, such a smoothing estimate was first shown to
Korteweg-de Vries equation
\[
\left\{
\begin{aligned}
&\partial_tu+\partial_x^3u+u\partial_xu=0,
\\
&u(0,x)=\varphi(x)\in L^2(\R),
\end{aligned}
\right.
\]
and Kato \cite{Ka2} proved that the solution $u=u(t,x)$ ($t,x\in\R$) satisfies
\[
\int^T_{-T}\int^R_{-R}\abs{\partial_xu(x,t)}^2\,dxdt
\leq c(T,R,\n{\varphi}_{L^2}).
\]
\par
Similar smoothing estimates have been observed for
generalised equations
\[
\left\{
\begin{aligned}
\p{i\partial_t+a(D_x)}\,u(t,x)&=0,\\
u(0,x)&=\varphi(x)\in L^2(\R^n),
\end{aligned}
\right.
\]
which come from equations of fundamental importance in mathematical physics
as their principal parts:
\begin{itemize}
 \item $a(\xi)=|\xi|^2$ $\cdots$ Schr\"odinger
\[
i\partial_tu-\Delta_xu=0
\]
 \item $a(\xi)=\sqrt{|\xi|^2+1}$ $\cdots$ Relativistic Schr\"odinger
\[
i\partial_tu+\sqrt{-\Delta_x+1}\,u=0
\]
 \item $a(\xi)=\xi^3$ ($n=1$) $\cdots$ Korteweg-de Vries (shallow water wave)
\[
\partial_tu+\partial_x^3u+u\partial_xu=0\
\]
 \item $a(\xi)=|\xi|\xi$ ($n=1$) $\cdots$ Benjamin-Ono (deep water wave)
\[
\partial_tu-\partial_x|D_x|u+u\partial_xu=0
\]
 \item $a(\xi)=\xi_1^2-\xi_2^2$ ($n=2$) $\cdots$ Davey-Stewartson
(shallow water wave of 2D)
\[
\begin{cases}
i\partial_tu-\partial_x^2u+\partial_y^2u=c_1|u|^2u+c_2u\partial_x v
\\
\partial_x^2 v-\partial_y^2 v=\partial_x|u|^2
\end{cases}
\]
\item
$a(\xi)=\xi_1^3+\xi_2^3,\, \xi_1^3+3\xi^2_2,\,
\xi_1^2+\xi_1\xi_2^2$ $\cdots$ Shrira (deep water wave of 2D)
\item $a(\xi)=\text{quadratic form}$ ($n\geq3$) $\cdots$ Zakharov-Schulman 
(interaction of sound wave and low amplitudes high frequency wave)
\end{itemize}

There has already been a lot of
literature on this subject from various points of view.
See, Ben-Artzi and Devinatz \cite{BD1, BD2},
Ben-Artzi and Klainerman \cite{BK}, Chihara \cite{Ch},
Hoshiro \cite{Ho1, Ho2},
Kato and Yajima \cite{KY},
Kenig, Ponce and Vega \cite{KPV1, KPV2, KPV3, KPV4, KPV5}, 
Linares and Ponce \cite{LP}, Simon \cite{Si}, 
Sugimoto \cite{Su1, Su2},
Walther \cite{Wa1, Wa2}, and many others.
We note that for a given operator $A$ the following are equivalent to each
other based on classical works by Agmon \cite{A} and Kato \cite{Ka1}:
\par
$\bullet$ Smoothing estimate
\[
\n{A e^{-it\Delta_x}\varphi(x)}
_{L^2\p{\R_t\times\R^n_x}}\le C\n{\varphi}_{L^2\p{\R^n_x}}
\quad\text{where $A=A(X,D_x)$},
\]
\par
$\bullet$ Restriction estimate
\[
\n{
\widehat{A^*f}_{|S^{n-1}_\rho}
}_{L^2\p{S^{n-1}_\rho}}
\le 
C\sqrt{\rho}
\n{f}_{L^2(\R^n)},
\quad\text{where $S^{n-1}_\rho=\b{\xi;\,|\xi|=\rho}$, ($\rho>0$)},
\]
\par
$\bullet$ Resolvent estimate
\[
 \sup_{\Im\zeta>0}\abs{\p{R(\zeta)A^*f,A^*f}} 
\leq
C \n{f}^2_{L^2(\R^n)},\quad
\text{where
$R(\zeta)=\p{-\triangle-\zeta}^{-1}$}.
\]
Most of the literature so far use the above equivalence
to show smoothing estimates
for dispersive equations
by showing restriction or resolvent estimates instead.

But here we develop a completely different strategy.
We investigate smoothing estimates by using methods
of comparison and
canonical transform which are quite efficient for this problem:
\par
\begin{enumerate}
 \item {\bf Comparison principle}
$\cdots$ comparison of symbols implies that of estimates,
 \item {\bf Canonical transform}
$\cdots$ transform an equation to another simple one.
\end{enumerate}
They work not only for all the dispersive equations
 (that is, the case $\nabla a\neq0$)
but also for some non-dispersive equations, and
induce smoothing estimates of an invariant form.
Smoothing estimates for inhomogeneous equations can be also
discussed by a similar treatment.
We will explain them in due order.

\section{Comparison principle}
Here we list theorems exemplifying the comparison principle,
which have been established in \cite[Section 2]{RS6}:
\begin{thm}[1D case]\label{1Dcomp}
Let $f,g\in C^1(\R)$ be real-valued and strictly
monotone. 
If $\sigma,\tau\in C^0(\R)$ satisfy
\[
\frac{|\sigma(\xi)|}{|f^\prime(\xi)|^{1/2}}
\leq A 
\frac{|\tau(\xi)|}{|g^\prime(\xi)|^{1/2}}
\]
then we have
\[
\|\sigma(D_x)e^{itf(D_x)}\varphi(x)\|_{L^2(\R_t)}
\leq A
\|\tau(D_x)e^{it g(D_x)}\varphi(x)\|_{L^2(\R_t)}
\]
for all $x\in\R$.
\end{thm}
\begin{thm}[2D case]\label{2Dcomp}
Let $f(\xi,\eta), g(\xi,\eta)\in C^1(\R^2)$ be real-valued and
strictly monotone in $\xi\in\R$ for each fixed $\eta\in\R$.
If $\sigma,\tau\in C^0(\R^2)$ satisfy
\[
\frac{|\sigma(\xi,\eta)|}
{\abs{f_\xi(\xi,\eta)}^{1/2}}
\leq A 
\frac{|\tau(\xi,\eta)|}
{\abs{g_\xi(\xi,\eta)}^{1/2}}
\]
then we have
\begin{multline*}
\n{\sigma(D_x,D_y)
e^{it f(D_x,D_y)}\varphi(x,y)}_{L^2(\R_t\times\R_{y})} \\ \leq A
\|\tau(D_x,D_y)e^{it g(D_x,D_y)}\varphi(x,y)\|_
{L^2(\R_t\times\R_{y})}
\end{multline*}
for all $x\in\R$.
\end{thm}
\begin{thm}[Radially Symmetric case]\label{Radcomp}
Let $f,g\in C^1(\R_+)$ be real-valued and strictly
monotone.
If $\sigma,\tau\in C^0(\R_+)$ satisfy
\[
\frac{|\sigma(\rho)|}{|f^\prime(\rho)|^{1/2}}
\leq A 
\frac{|\tau(\rho)|}{|g^\prime(\rho)|^{1/2}}
\]
then we have
\[
\|\sigma(|D_x|)e^{it f(|D_x|)}
\varphi(x)\|_{L^2(\R_t)}
\leq A
\|\tau(|D_x|)e^{it g(|D_x|)}\varphi(x)\|_{L^2(\R_t)}
\]
for all $x\in\R^n$.
\end{thm}

\section{Canonical transforms}
\label{SECTION:canonical}
Next we will review the idea of canonical transforms
discussed in \cite[Section 4]{RS6}.
It is based on the so-called Egorov's theorem.

Let $\psi:\Gamma\to\widetilde{\Gamma}$
be a $C^\infty$-diffeomorphism between open sets
$\Gamma\subset\R^n$ and $\widetilde{\Gamma}\subset\R^n$.
We always assume that
\[
C^{-1}\leq\abs{\det \partial\psi(\xi)}\leq C\quad(\xi\in\Gamma),
\]
for some $C>0$.
We set formally
\[
I_\psi u(x)
=\FT^{-1}\left[\FT u\p{\psi(\xi)}\right](x)
=(2\pi)^{-n}\int_{\R^n}\int_{\R^n}
   e^{i(x\cdot\xi-y\cdot\psi(\xi))}u(y)\, dyd\xi.
\]
The operators $I_\psi$ can be justified
by using cut-off functions
$\gamma\in C^\infty(\Gamma)$ and
$\widetilde{\gamma}=\gamma\circ\psi^{-1}\in C^\infty(\widetilde{\Gamma})$
which satisfy $\supp\gamma\subset\Gamma$,
$\supp\widetilde{\gamma}\subset\widetilde{\Gamma}$.
We set
\begin{equation}\label{DefI0}
\begin{aligned}
I_{\psi,\gamma} u(x)
&=\FT^{-1}\left[\gamma(\xi)\FT u\p{\psi(\xi)}\right](x)
\\
&=(2\pi)^{-n}\int_{\R^n}\int_{\Gamma}
 e^{i(x\cdot\xi-y\cdot\psi(\xi))}\gamma(\xi)u(y) dyd\xi.
\end{aligned}
\end{equation}
In the case that $\Gamma$, $\widetilde{\Gamma}\subset\R^n\setminus0$
are open cones,
we may consider the homogeneous functions $\psi$ and $\gamma$ which satisfy
$\supp\gamma\cap \Sph^{n-1}\subset\Gamma\cap \Sph^{n-1}$ and
$\supp\widetilde{\gamma}\cap 
\Sph^{n-1}\subset\widetilde{\Gamma}\cap \Sph^{n-1}$,
where $\Sph^{n-1}=\b{\xi\in\Rn: |\xi|=1}$.
Then we have the expressions for compositions
\[
I_{\psi,\gamma}=\gamma(D_x)\cdot I_\psi=I_\psi\cdot\widetilde{\gamma}(D_x)
\]
and also the formula
\begin{equation}\label{eq:cnon}
I_{\psi,\gamma}\cdot\sigma(D_x)
=\p{\sigma\circ\psi}(D_x)\cdot I_{\psi,\gamma}.
\end{equation}
\par
We also introduce the weighted $L^2$-spaces.
For a weight function $w(x)$, let $L^2_{w}(\R^n;w)$ be
the set of measurable functions $f:\Rn\to\C$ 
such that the norm
\[
\n{f}_{L^2(\R^n;w)}
=\p{\int_{\R^n}\abs{w(x) f(x)}^2\,dx}^{1/2}
\]
is finite.
Then, on account of the relations \eqref{eq:cnon},
we obtain the following fundamental theorem (\cite[Theorem 4.1]{RS6}):
\begin{thm}\label{Th:reduction}
Assume that the operator $I_{\psi,\gamma}$ defined by \eqref{DefI0}
is $L^2(\R^n;w)$--bounded.
Suppose that we have the estimate
\[
\n{w(x)\rho(D_x)e^{it\sigma(D_x)}\varphi(x)}_{L^2\p{\R_t\times\R^n_x}}
\leq C\n{\varphi}_{L^2\p{\R^n_x}}
\]
for all $\varphi$ such
that $\supp\widehat{\varphi}\subset\supp\widetilde\gamma$, where
$\widetilde\gamma=\gamma\circ\psi^{-1}$.
Assume also that the function
\[
q(\xi)=\frac{\gamma\cdot\zeta}{\rho\circ \psi}(\xi)
\]
is bounded.
Then we have
\[
\n{w(x)\zeta(D_x)e^{ita(D_x)}\varphi(x)}_{L^2\p{\R_t\times\R^n_x}}
\leq C\n{\varphi}_{L^2\p{\R^n_x}}
\]
for all $\varphi$
such that $\supp\widehat{\varphi}\subset\supp\gamma$,
where $a(\xi)=(\sigma\circ\psi)(\xi)$.
\end{thm}
Note that $e^{ita(D_x)}\varphi(x)$ and $e^{it\sigma(D_x)}\varphi(x)$ are
solutions to
\[
\left\{
\begin{aligned}
\p{i\partial_t+a(D_x)}\,u(t,x)&=0,\\
u(0,x)&=\varphi(x),
\end{aligned}
\right.
\quad {\rm and}\quad
\left\{
\begin{aligned}
\p{i\partial_t+\sigma(D_x)}\,v(t,x)&=0,\\
v(0,x)&=g(x),
\end{aligned}
\right.
\]
respectively.
Theorem \ref{Th:reduction} means that smoothing estimates for the equation with
$\sigma(D_x)$ implies those with $a(D_x)$ if the canonical transformations
which
relate them are bounded on weighted $L^2$-spaces.

As for the $L^2(\R^n;w)$--boundedness of the operator $I_{\psi,\gamma}$,
we have criteria for some special weight functions.
For $\ka\in\R$, let $L^2_\ka(\R^n)$ be
the set of measurable functions $f$ such that the norm
\[
\n{f}_{L^2_\ka(\R^n)}
=\p{\int_{\R^n}\abs{\langle x\rangle^\ka f(x)}^2\,dx}^{1/2}
\]
is finite.
Then we have the following mapping properties
(\cite[Theorems 4.2, 4.3]{RS6}).
\begin{thm}\label{Th:L'2k}
Let $\Gamma$, $\widetilde{\Gamma}\subset\R^n\setminus0$ be open cones.
Suppose $|\ka|< n/2$.
Assume $\psi(\lambda\xi)=\lambda\psi(\xi)$,
$\gamma(\lambda\xi)=\gamma(\xi)$ for all $\lambda>0$ and $\xi\in\Gamma$.
Then the operator $I_{\psi,\gamma}$
defined by \eqref{DefI0} is $L^2_{\ka}(\R^n)$--bounded.
\end{thm}
\begin{thm}\label{Th:L2k}
Suppose $\ka\in\R$.
Assume that all the derivatives of entries of the 
$n\times n$ matrix
$\partial\psi$ and those of $\gamma$ are bounded.
Then the operator $I_{\psi,\gamma}$ 
defined by
\eqref{DefI0} are $L^2_{\ka}(\R^n)$--bounded.
\end{thm}

\section{Smoothing estimates for dispersive equations}
\par
We consider smoothing estimates for solutions
$u(t,x)=e^{ita(D_x)}\varphi(x)$
to general equations
\[
\left\{
\begin{aligned}
\p{i\partial_t+a(D_x)}\,u(t,x)&=0,\\
u(0,x)&=\varphi(x)\in L^2(\R^n).
\end{aligned}
\right.
\]
Let $a_m(\xi)$ be the principal term of $a(\xi)$ satisfying
\[
 a_m(\xi)\in C^\infty(\R^n\setminus0),\quad\text{real-valued},
\quad a_m(\lambda\xi)=\lambda^ma_m(\xi)\,\,\, (\lambda>0,\xi\neq0).
\]
We assume that $a(\xi)$ is {\it dispersive} in the following sense:
\begin{equation}\tag{{\bf H}}
a(\xi)=a_m(\xi),\qquad\nabla a_m(\xi)\neq0 \quad(\xi\in\R^n\setminus0),
\end{equation}
or, otherwise, we assume
\begin{equation}\tag{{\bf L}}
\begin{aligned}
&a(\xi)\in C^\infty(\R^n),\qquad \nabla a(\xi)\neq0 \quad(\xi\in\R^n),
\quad\nabla a_m(\xi)\neq0 \quad (\xi\in\R^n\setminus0), 
\\
&|\partial^\alpha\p{a(\xi)-a_m(\xi)}|\leq C_\alpha\abs{\xi}^{m-1-|\alpha|}
\quad\text{for all multi-indices $\alpha$ and all $|\xi|\geq1$}.
\end{aligned}
\end{equation}
\begin{example}
$a(\xi)=\xi_1^3+\cdots+\xi_n^3+\xi_1$ satisfies {\rm (L)}.
\end{example}
\par
The dispersiveness means that the classical orbit, that is, the solution of
the Hamilton-Jacobi equations 
\[
\left\{
\begin{aligned}
\dot{x}(t)&=\p{\nabla a}(\xi(t)), \quad\dot{\xi}(t)=0,
\\
x(0)&=0,\quad \xi(0)=k,
\end{aligned}
\right.
\]
does not stop, and the singularity of $u(t,x)=e^{ita(D_x)}\varphi(x)$
travels to infinity along this orbit.
Hence we can expect the smoothing, and indeed we have
the following result (\cite[Theorem 5.1, Corollary 5.5]{RS6}):
\begin{thm}\label{Th2.1}
Assume {\rm (H)} or {\rm (L)}.
Suppose $m\geq1$ and $s>1/2$. Then we have
\[
\n{\langle x\rangle^{-s}|D_x|^{(m-1)/2} e^{ita(D_x)}\varphi(x)}
_{L^2\p{\R_t\times\R^n_x}}
\le C\n{\varphi}_{L^2\p{\R^n}}.
\]
\end{thm}
\begin{remark}\label{Rem2.1}
Theorem \ref{Th2.1} with polynomials $a(\xi)$
follows immediately from a sharp version of 
local smoothing estimate proved by Kenig, Ponce and Vega
\cite[Theorem 4.1]{KPV1},
and any polynomial $a(\xi)$ which satisfies the estimate in Theorem \ref{Th2.1}
has to be dispersive, that is 
$\nabla a_m(\xi)\neq0$ $(\xi\neq0)$
(see Hoshiro \cite{Ho2}).
Theorem \ref{Th2.1} with $a(\xi)=|\xi|^2$ and $n\geq 3$
was also stated by Ben-Artzi and Klainerman \cite{BK},
and with the case (H) and $m>1$ by Chihara \cite{Ch} in different contexts.
\end{remark}

\section{Proof by new methods}
We explain how to prove Theorem \ref{Th2.1} under the condition (H)
by our new method.
The main strategy is that we obtain estimates for low dimensional model cases
from some trivial estimate by the comparison principle,
and reduce general case to such model cases by the method
of canonical transforms.

\subsection{Low dimensional model estimates}
By the comparison principle, we can show the equivalence
of low dimensional estimates of various type.
In the 1D case, we have (for $l,m>0$)
\begin{equation}\label{dim1ex}
\sqrt m\,
\n{|D_x|^{(m-1)/2}e^{it|D_x|^{m}}
\varphi(x)}_{L^2(\R_t)}
=
\sqrt l\,
\n{|D_x|^{(l-1)/2}e^{it|D_x|^{l}}
\varphi(x)}_{L^2(\R_t)}
\end{equation}
for all $x\in\R$.
Here $\supp\widehat{\varphi}\subset [0,+\infty)$ or $(-\infty,0]$.

In the 2D case, we have (for $l,m>0$)
\begin{multline}\label{dim2ex}
\n{|D_y|^{(m-1)/2}
e^{itD_x|D_y|^{m-1}}\varphi(x,y)}_{L^2(\R_t\times\R_y)}
\\=
\n{|D_y|^{(l-1)/2}e^{it{D_x|D_y|^{l-1}}\varphi(x,y)}}
_{L^2(\R_t\times\R_y)}
\end{multline}
for all $x\in\R$.
On the other hand, in 1D case, we have trivially
\begin{equation}\label{core}
\n{e^{itD_x}\varphi(x)}_{L^2(\R_t)}
=
\n{\varphi(x+t)}_{L^2\p{\R_{x}}}
=\n{\varphi}_{L^2\p{\R_{x}}}
\end{equation}
for all $x\in\R$.
Using the equality \eqref{core},
the right hand sides of \eqref{dim1ex} and \eqref{dim2ex} with $l=1$ can be
estimated, and we have for all $x\in\R$:
\par
$\bullet$ (1D Case)
\[
\n{|D_x|^{(m-1)/2}e^{it|D_x|^m}\varphi(x)}_{L^2(\R_t)}
\leq C\n{\varphi}_{L^2(\R_x)},
\]
\par
$\bullet$ (2D Case)
\[
\n{|D_y|^{(m-1)/2}e^{itD_x|D_y|^{m-1}}\varphi(x,y)}_{L^2(\R_t\times\R_y)}
\leq C\n{\varphi}_{L^2\p{\R^2_{x,y}}}.
\]
\begin{remark}
In the case $m=2$,
these estimates were proved by
Kenig, Ponce \& Vega \cite{KPV1} (1D case)
and Linares \& Ponce \cite{LP} (2D case).
\end{remark}
\par
The following is a straightforward consequence from these estimates:
\begin{proposition}\label{Prop5.1}
Suppose $m>0$ and $s>1/2$.
Then for $n\geq1$ we have
\[
\n{\jp{x}^{-s}|D_n|^{(m-1)/2}
e^{it|D_n|^m}\varphi(x)}_{L^2(\R_t\times\R^n_x)}
\leq
 C\n{\varphi}_{L^2(\R_x^n)}
\]
and for $n\geq2$ we have
\[
\n{\jp{x}^{-s}|D_n|^{(m-1)/2}e^{itD_1|D_n|^{m-1}}\varphi(x)}
_{L^2(\R_t\times\R^n_x)}
\leq
 C\n{\varphi}_{L^2(\R_x^n)},
\]
where $D_x=(D_1,\ldots,D_n)$.
\end{proposition}

\subsection{Reduction to model estimates}
On account of the method of canonical transform (Theorem \ref{Th:reduction}),
smoothing estimates for dispersive equations (Theorem \ref{Th2.1})
can be reduced
to low dimensional model estimates (Proposition \ref{Prop5.1})
by the canonical transformation
if we find a homogeneous change of variable $\psi$
such that
\[
a(\xi)=\p{\sigma\circ\psi}(\xi),
\quad
\sigma(D)=|D_n|^m\quad\text{or}\quad \sigma(D)=D_1|D_n|^{m-1}.
\]
We show how to select such $\psi$ under the assumption (H).
The argument for the case (L) is similar.
By microlocalisation and rotation,
we may assume that the initial data $\varphi$
satisfies $\supp \hat{\varphi}\subset\Gamma$, where
$\Gamma\subset\R^n\setminus0$ is a sufficiently small conic neighbourhood of
$e_n=(0,\ldots0,1)$.
Furthermore, we have Euler's identity
\[
a(\xi)=a_m(\xi)=\frac1m\xi\cdot\nabla a(\xi),
\]
and the dispersiveness $\nabla a(e_n)\neq0$
implies the following two cases:
\begin{description}
\item[(I)] 
\hspace{.2mm}
$\partial_n a(e_n)\neq0$ $\cdots$ (elliptic).
By Euler's identity,
we have $a(e_n)\neq0$.
Hence, in this case, we may assume
$a(\xi)>0$ $(\xi\in\Gamma)$, $\partial_n a(e_n)\neq0$.
\item[(II)]
$\partial_n a(e_n)=0$ $\cdots$ (non-elliptic).
By assumption $\nabla a(e_n)\neq0$, there exits $j\neq n$ such that
$\partial_j a(e_n)\neq0$.
Hence, in this case, we may assume
$\partial_1a(e_n)\neq0$.
\end{description}
In the elliptic case (I), we take
\[
\sigma(\eta)=|\eta_n|^m,\quad
\psi(\xi)=(\xi_1,\ldots,\xi_{n-1},a(\xi)^{1/m}).
\]
Then we have $a(\xi)=\p{\sigma\circ\psi}(\xi)$, and $\psi$ is surely a change of 
variables on $\Gamma$ since
\[
\det\partial\psi(e_n)
=
\begin{vmatrix}
E_{n-1}&0
\\
*&\frac1m a(e_n)^{1/m-1}\partial_n a(e_n)
\end{vmatrix}
\neq0
\]
where $E_{n-1}$ is the identity matrix.
In the non- elliptic case (II), we take
\[
\sigma(\eta)=\eta_1|\eta_n|^{m-1},\quad
\psi(\xi)=\p{\frac{a(\xi)}{|\xi_n|^{m-1}},\xi_2,\ldots,\xi_n}.
\]
Then we have  again $a(\xi)=\p{\sigma\circ\psi}(\xi)$ and
\[
\det\partial\psi(e_n)
=
\begin{vmatrix}
\partial_1 a(e_n)&*
\\
0&E_{n-1}
\end{vmatrix}
\neq0.
\]
Thus, we successfully showed Theorem \ref{Th2.1}
in both cases.

\section{Non-dispersive case}
Now we consider what happens
if the equation does not satisfy the dispersiveness assumption
$\nabla a(\xi)\neq0$ $(\xi\in\R^n)$.
All the precise results and arguments in this section are to appear in
our forthcoming paper \cite{RS2}.

Although we cannot have smoothing estimates (see Remark \ref{Rem2.1}),
such case appears naturally in physics.
For example, let us consider a coupled system of Schr\"odinger equations
\[
i\partial_t v  =\Delta_x v+b(D_x)w, \quad
i\partial_t w  =\Delta_x w+c(D_x)v, 
\]
which represents a linearised model of wave packets with two modes.
Assume that this system is diagonalised and regard it
as a single equations for the eigenvalues:
\[
a(\xi)=-|\xi|^2\pm\sqrt{b(\xi)c(\xi)}.
\]
Then there could exist points $\xi$ such that $\nabla a(\xi)=0$
because of the lower order terms $b(\xi)$, $c(\xi)$.
Another interesting examples are Shrira equations, in which case:
\[
a(\xi)=\xi_1^3+\xi_2^3,\quad \xi_1^3+3\xi^2_2,\quad
\xi_1^2+\xi_1\xi_2^2.
\]
Although $a(\xi)=\xi_1^3+\xi_2^3$ satisfies assumption (H),
$a(\xi)=\xi_1^3+3\xi^2_2$ and $a(\xi)=\xi_1^2+\xi_1\xi_2^2$ do not
satisfy assumption (L) because $\nabla a(0)=0$.

\par
We suggest an estimate which we expect to hold for non-dispersive equations:
\begin{equation}\label{invariant}
\n{\jp{x}^{-s}|\nabla a(D_x)|^{1/2}
e^{it a(D_x)}\varphi(x)}_{L^2\p{\R_t\times\R^n_x}}
\leq C\n{\varphi}_{L^2\p{\R^n_x}} \quad {\rm(}s>1/2{\rm)}
\end{equation}
and let us call it {\it invariant estimate}.
This estimate has a number of advantages:
\begin{itemize}
\item in the dispersive case $\nabla a(\xi)\neq0$,
it is equivalent to Theorem \ref{Th2.1};

\item it is invariant under canonical transformations
for the operator $a(D_x)$;

\item it does continue to hold for a variety of non-dispersive
operators $a(D_x)$, where $\nabla a(\xi)$ may become zero on some set
and when the usual estimate fails;

\item it does take into account zeros of the gradient
$\nabla a(\xi)$, which is also responsible for the interface
between dispersive and non-dispersive zone (e.g. how quickly
the gradient vanishes);
\end{itemize}

\subsection{Secondary comparison}
\par
By using comparison principle again to the smoothing estimates
obtained from the comparison principle, we can have new estimates.
This is a powerful tool to induce the invariant estimates \eqref{invariant} for
non-dispersive equations.
For example, we have just obtained the estimate
\[
\n{\langle x\rangle^{-s}|D_x|^{(m-1)/2} e^{it|D_x|^m}\varphi}
_{L^2\p{\R_t\times\R^n_x}}
\le C\n{\varphi}_{L^2\p{\R^n_x}}
\]
(Theorem \ref{Th2.1} with $a(\xi)=|\xi|^m$) from comparison principle
and canonical transformation.
If we set
$g(\rho)=\rho^m$, $\tau(\rho)=\rho^{(m-1)/2}$, then we have
$|\tau(\rho)|/|g'(\rho)|^{1/2}=1/\sqrt{m}$.
Hence by the comparison result again for the radially symmetric case
(Theorem \ref{Radcomp}), we have
\begin{thm}
Suppose $s>1/2$.
Let $f\in C^1(\R_+)$ be real-valued and strictly monotone.
If $\sigma\in C^0(\R_+)$ satisfy
\[
|\sigma(\rho)|
\leq A|f^\prime(\rho)|^{1/2},
\]
then we have
\[
\|\jp{x}^{-s}\sigma(|D_x|)e^{it f(|D_x|)}\varphi(x)\|
_{L^2\p{\R_t\times\R^n_x}}
\leq C\n{\varphi}_{L^2\p{\R^n_x}}.
\]
\end{thm}
From this secondary comparison, we obtain immediately the following invariant
estimate since a radial function $a(\xi)=f(|\xi|)$ always satisfies
$|\nabla a(\xi)|=|f'(|\xi|)|$.
\begin{thm}\label{SymInv}
Suppose $s>1/2$.
Let $a(\xi)=f(|\xi|)$ and $f\in C^\infty(\R_+)$ be real-valued.
Then we have
\[
\n{\jp{x}^{-s}|\nabla a(D_x)|^{1/2}
e^{it a(D_x)}\varphi(x)}_{L^2\p{\R_t\times\R^n_x}}
\leq C\n{\varphi}_{L^2\p{\R^n_x}}.
\]
\end{thm}
\begin{example}
$a(\xi)=(|\xi|^2-1)^2$ is non-dispersive
because $$\nabla a(\xi)=4(|\xi|^2-1)\xi=0$$ if
$|\xi|=0,1$.
But we have the invariant estimate by Theorem \ref{SymInv}.
\end{example}
\par
For the non-radially symmetric case, we compare again to the low dimensional model
estimates (Proposition \ref{Prop5.1})
and obtain
\begin{thm}[1D secondary comparison]\label{1Dseccomp}
Suppose $s>1/2$.
Let $f\in C^1(\R)$ be real-valued and strictly
monotone.
If $\sigma\in C^0(\R)$ satisfies
\[
|\sigma(\xi)|
\leq A |f^\prime(\xi)|^{1/2},
\]
then we have
\[
\|\jp{x}^{-s}\sigma(D_x)e^{it f(D_x)}\varphi(x)\|_{L^2(\R_t\times\R_x)}
\leq AC
\|\varphi(x)\|_{L^2(\R_x)}.
\]
\end{thm}
\begin{thm}[2D secondary comparison]\label{2Dseccomp}
Suppose $s>1/2$.
Let $f\in C^1(\R^2)$ be real-valued and $f(\xi,\eta)$ be
strictly monotone in $\xi\in\R$ for every fixed $\eta\in\R$.
If $\sigma\in C^0(\R^2)$ satisfies
\[
|\sigma(\xi,\eta)|
\leq A \abs{\partial f/\partial \xi(\xi,\eta)}^{1/2},
\]
then we have
\[
\n{\jp{x}^{-s}\sigma(D_x,D_y)
e^{it f(D_x,D_y)}\varphi(x,y)}_{L^2(\R_t\times\R_{x,y}^2)}  \leq AC
\|\varphi(x,y)\|_{L^2(\R_{x,y}^2)}.
\]
\end{thm}
\begin{example}
By using secondary comparison for non-radially symmetric case,
we have invariant estimates for Shrira equations.
In fact, for $a(\xi)=\xi_1^3+3\xi_2^2$,
we have by 1D secondary comparison (Theorem \ref{1Dseccomp})
\begin{align*}
&\n{\jp{x_1}^{-s}|D_1|
e^{it D_1^3}\varphi(x)}_{L^2(\R_t\times\R_x^2)}
\leq C\n{\varphi}_{L^2(\R_x^2)}, \\
&\n{\jp{x_2}^{-s}|D_2|^{1/2}
e^{it 3D_2^2}\varphi(x)}_{L^2(\R_t\times\R_x^2)}
\leq C\n{\varphi}_{L^2(\R_x^2)},
\end{align*}
for $s>1/2$.
Hence by $\jp{x}^{-s}\leq\jp{x_k}^{-s}$ ($k=1,2$)
we have
\[
\n{\jp{x}^{-s}\p{|D_1|+|D_2|^{1/2}}
e^{it a(D_x)}\varphi(x)}_{L^2(\R_t\times\R_x^2)}
\leq C\n{\varphi}_{L^2(\R_x^2)}
\]
and hence we have
\[
\n{\jp{x}^{-s}|\nabla a(D_x)|^{1/2}
e^{it a(D_x)}\varphi(x)}_{L^2(\R_t\times\R_x^2)}
\leq C\n{\varphi}_{L^2(\R_x^2)}.
\]
For $a(\xi)=\xi_1^2+\xi_1\xi_2^2$, we have
by 2D secondary comparison (Theorem \ref{2Dseccomp})
\begin{align*}
&\n{\jp{x_1}^{-s}|2D_1+D_2^2|^{1/2}
e^{it a(D_1,\,D_2)}\varphi(x)}_{L^2(\R_t\times\R_x^2)}
\leq C\n{\varphi}_{L^2(\R_x^2)},
\\
&\n{\jp{x_2}^{-s}|D_1D_2|^{1/2}
e^{it a(D_1,\,D_2)}\varphi(x)}_{L^2(\R_t\times\R_x^2)}
\leq C\n{\varphi}_{L^2(\R_x^2)},
\end{align*}
for $s>1/2$, hence we have similarly
\[
\n{\jp{x}^{-s}|\nabla a(D_x)|^{1/2}
e^{it a(D_x)}\varphi(x)}_{L^2(\R_t\times\R_x^2)}
\leq C\n{\varphi}_{L^2(\R_x^2)}.
\]
\end{example}

\subsection{Non-dispersive case controlled by Hessian}
We will show that in the non-dispersive situation
the rank of $\nabla^2 a(\xi)$ still has a responsibility for
smoothing properties.

First let us consider the case when dispersiveness (L) is true only
for large $\xi$:
\begin{equation}\tag{{\bf L'}}
\begin{aligned}
&|\nabla a(\xi)|\geq C\jp{\xi}^{m-1}\quad(|\xi|\gg1),
\\
&|\partial^\alpha\p{a(\xi)-a_m(\xi)}|\leq C\jp{\xi}^{m-1-|\alpha|}
\quad(|\xi|\gg1).
\end{aligned}
\end{equation}
\begin{thm}\label{isolated-critical}\label{nondeg}
Suppose $n\geq1$, $m\geq1$, and $s>1/2$.
Let $a\in C^\infty(\Rn)$ be real-valued and assume that it
has finitely many critical points.
Assume {\rm{(L')}} and
\[
\nabla a(\xi)=0\, \Rightarrow\, \det\nabla^2 a(\xi)\not=0. 
\]
Then we have
\[
\n{\jp{x}^{-s}|\nabla a(D_x)|^{1/2}
e^{it a(D_x)}\varphi(x)}_{L^2\p{\R_t\times\R^n_x}}
\leq C\n{\varphi}_{L^2\p{\R^n_x}}.
\]
\end{thm}
\begin{example}
$a(\xi)=\xi_1^4+\cdots+\xi_n^4+|\xi|^2$ satisfies assumptions
in Theorem \ref{isolated-critical}.
\end{example}
We outline the proof of Theorem \ref{nondeg}.
For the region where $\nabla a(\xi)\neq0$, we can use 
a smoothing estimate for dispersive equations.
Near the points $\xi$ where $\nabla a(\xi)=0$,
there exists a change of variable $\psi$ by Morse's lemma
such that
$a(\xi)=(\sigma\circ\psi)(\xi)$ where
$\sigma(\eta)$ is a non-degenerate quadratic form,
and satisfies dispersiveness (H).
Hence the estimate can be reduced to the dispersive case
by the method of canonical transformation.
\par
Next we consider the case when $a(\xi)$ is homogeneous (of oder $m$).
Then, by Euler's identity, we have
\[
\nabla a(\xi)=\frac1{m-1}\xi\nabla^2 a(\xi) \quad(\xi\neq0),
\]
hence
\[
\nabla a(\xi)=0 \,\Rightarrow\, \det\nabla^2a(\xi)= 0\quad(\xi\neq0).
\]
Therefore assumption in Theorem \ref{isolated-critical}
does not make any sense in this case,
but we can have the following result if we use the idea of canonical
transform again:
\begin{thm}\label{th1}
Suppose $n\geq2$ and $s>1/2$.
Let $a\in C^\infty(\Rn\setminus0)$ be real-valued and satisfy
$a(\lambda\xi)=\lambda^2 a(\xi)$ $(\lambda>0,\,\xi\neq0)$.
Assume that
\[
\nabla a(\xi)=0 \,\Rightarrow\, \rank\nabla^2a(\xi)= n-1\quad (\xi\neq0).
\]
Then we have
\[
\n{\jp{x}^{-s}|\nabla a(D_x)|^{1/2}
e^{it a(D_x)}\varphi(x)}_{L^2\p{\R_t\times\R^n_x}}
\leq C\n{\varphi}_{L^2\p{\R^n_x}}.
\]
\end{thm}
\begin{example}
$\displaystyle{
a(\xi)=\frac{\xi_1^2\xi_2^2}{\xi_1^2+\xi_2^2}
+\xi^2_3+\cdots+\xi_n^2}$
satisfies the assumptions in Theorem \ref{th1}.
In the case $n=2$, this is an illustration of a smoothing estimate for the 
Cauchy problem for an equation like
$$
i\partial_t u+D_1^2D_2^2 \Delta^{-1}u=0
$$
which is regarded as a mixture of Davey-Stewartson
and Benjamin-Ono type equations.
\end{example}

\section{Concluding remarks}
\subsection{Summary}
Finally we summarise what is explained in this article in a diagram below.
It is remarkable that all the results of smoothing estimates so far
is derived from just the
translation invariance of Lebesgue measure:
\par
\bigskip
\par
$\bullet$ Trivial estimate
$\n{\varphi(x+t)}_{L^2(\R_t)}=\n{\varphi}_{L^2\p{\R_{x}}}$
\medskip
\[
\qquad\qquad\Downarrow\qquad\text{(comparison principle)}
\]
\par
$\bullet$ Low dimensional model estimates (Proposition \ref{Prop5.1})
\medskip
\[
\qquad\qquad\Downarrow\qquad\text{(canonical transform)}
\]
\par
$\bullet$ Smoothing estimates for dispersive equations (Theorem \ref{Th2.1})
\[
\qquad\qquad\Downarrow\qquad\text{(secondary comparison \& canonical transform)}
\]
\par
$\bullet$ Invariant estimates for non-dispersive equations
at least for
\par
\medskip
\par\hspace{5mm}
$*$ radially symmetric $a(\xi)=f(|\xi|)$, $f\in C^1(\R_+)$,
\par
\medskip
\par\hspace{5mm}
$*$ Shrira equation $a(\xi)=\xi_1^3+3\xi^2_2,\,\,\xi_1^2+\xi_1\xi_2^2$,
\par
\medskip
\par\hspace{5mm}
$*$ non-dispersive $a(\xi)$ controlled by its Hessian.

\subsection{Smoothing estimates for inhomogeneous equations}

We finish this article by
mentioning some results for inhomogeneous equations.
Let us consider the solution
\[
u(t,x)
=-i\int^t_0e^{i(t-\tau)a(D_x)}f(\tau,x)\,d\tau
\]
to the equation
\[
\left\{
\begin{aligned}
\p{i\partial_t+a(D_x)}\,u(t,x)&=f(t,x)\quad\text{in $\R_t\times\R^n_x$},\\
u(0,x)&=0\quad\text{in $\R^n_x$}.
\end{aligned}
\right.
\]
Although smoothing estimates for such equation
are necessary for nonlinear applications
(see \cite{RScmp} for example), there are considerably
less results on this topic available in the literature.
But the method of canonical transform also works to this problem,
and we will list here some recent achievement 
given in our forthcoming paper \cite{RS7}.
The following result is a counter part of Theorem \ref{Th2.1}.
Especially, this kind of time-global estimate 
for the operator $a(D_x)$ with lower order terms are the benefit
of our new method:
\begin{thm}\label{Cor:inhom2}
Assume {\rm (H)} or {\rm (L)}.
Suppose $n\geq2$, $m\geq1$, and $s>1/2$.
Then we have
\[
\n{\jp{x}^{-s}|D_x|^{m-1}\int^t_0e^{i(t-\tau)a(D_x)}f(\tau,x)\,d\tau}_
{L^2(\R_t\times\R^n_x)}
\leq
 C\n{\jp{x}^s f(t,x)}_{L^2(\R_t\times\R^n_x)}.
\]
\end{thm}
The proof of Theorem \ref{Cor:inhom2} is carried out
by reducing it to model estimates below via canonical transform:
\begin{proposition}\label{Prop:inhom}
Suppose $n=1$ and $m>0$.
Let $a(\xi)\in C^\infty\p{\R\setminus0}$ be a real-valued
function which satisfies $a(\lambda\xi)=\lambda^m a(\xi)$ for
all $\lambda>0$ and $\xi\neq0$.
Then we have
\[
\n{a'(D_x)
\int^t_0e^{i(t-\tau)a(D_x)}f(\tau,x)\,d\tau}_{L^2(\R_t)}
\leq C\int_\R \n{f(t,x)}_{L^2(\R_t)}\,dx
\]
for all $x\in\R$.
Suppose $n=2$ and $m>0$.
Then we have
\[
\n{|D_x|^{m-1}\int^t_0e^{i(t-\tau)|D_x|^{m-1}D_y}
f(\tau,x,y)\,d\tau}_{L^2(\R_t\times\R_x)}
\\
\leq
C\int_\R \n{f(t,x,y)}_{L^2(\R_t\times\R_x)}\,dy
\]
for all $y\in\R$.
\end{proposition}

\begin{remark}
Proposition \ref{Prop:inhom} with the case $n=1$
is a unification of the results
by Kenig, Ponce and Vega who treated the
cases $a(\xi)=\xi^2$ (\cite[p.258]{KPV3}), $a(\xi)=|\xi|\xi$
(\cite[p.160]{KPV4}),
and $a(\xi)=\xi^3$ (\cite[p.533]{KPV2}).
\end{remark}
Since we unfortunately do not know the comparison principle
for inhomogeneous equations, we gave a direct proof to
Proposition \ref{Prop:inhom} in \cite{RS7}.

\end{document}